\documentclass{amsart}
\usepackage{amsmath}
\usepackage{amsfonts}

\newtheorem{theorem}{Theorem}
\theoremstyle{plain}

\newtheorem{corollary}{Corollary}

\numberwithin{equation}{section}
\input{tcilatex}

\begin{document}
\title[The $\left( h,q\right) $\textbf{\ numbers and polynomials with weight 
}$\alpha $]{\textbf{Some new identities on the (}$h,q$\textbf{)-Genocchi
numbers and polynomials with weight }$\alpha $}
\author[\textbf{S. Araci}]{\textbf{S. Araci}}
\address{\textbf{University of Gaziantep, Faculty of Science and Arts,
Department of Mathematics, 27310 Gaziantep, TURKEY}}
\email{\textbf{mtsrkn@hotmail.com}}
\author[\textbf{E. Cetin}]{\textbf{E. Cetin}}
\address{\textbf{Uludag University, Faculty of Arts and Science, Department
of Mathematics, Bursa, Turkey}}
\email{\textbf{elifc2@hotmail.com}}
\author[\textbf{M. Acikgoz}]{\textbf{M. Acikgoz}}
\address{\textbf{University of Gaziantep, Faculty of Science and Arts,
Department of Mathematics, 27310 Gaziantep, TURKEY}}
\email{\textbf{acikgoz@gantep.edu.tr}}
\author[\textbf{I. N. Cangul}]{\textbf{I. N. Cangul}}
\address{\textbf{Uludag University, Faculty of Arts and Science, Department
of Mathematics, Bursa, Turkey}}
\email{\textbf{ncangul@gmail.com}}

\begin{abstract}
We give some new identities for ($h,q$)-Genocchi numbers and polynomials by
means of the fermionic $p$-adic $q$-integral on $%
\mathbb{Z}
_{p}$ and the weighted $q$-Bernstein polynomials.

\vspace{2mm}\noindent \textsc{2000 Mathematics Subject Classification.}
05A10, 11B65, 28B99, 11B68, 11B73.

\vspace{2mm}

\end{abstract}

\maketitle

\section{\textbf{Introduction and Notations}}

Let $p$ be a fixed odd prime number. Throughout this paper we use the
following notations. By $%
\mathbb{Z}
_{p}$ we denote the ring of $p$-adic rational integers, $%
\mathbb{Q}
$ denotes the field of rational numbers, $%
\mathbb{Q}
_{p}$ denotes the field of $p$-adic rational numbers, and $%
\mathbb{C}
_{p}$ denotes the completion of algebraic closure of $%
\mathbb{Q}
_{p}$. Let $%
\mathbb{N}
$ be the set of natural numbers and $%
\mathbb{N}
^{\ast }=%
\mathbb{N}
\cup \left\{ 0\right\} $. The $p$-adic absolute value is defined by 
\begin{equation*}
\left\vert p\right\vert _{p}=\frac{1}{p}\text{.}
\end{equation*}%
In this paper, we assume $\left\vert q-1\right\vert _{p}<1$ as an
indeterminate. Let $UD\left( 
\mathbb{Z}
_{p}\right) $ be the space of uniformly differentiable functions on $%
\mathbb{Z}
_{p}$. For $f\in UD\left( 
\mathbb{Z}
_{p}\right) $, the fermionic $p$-adic $q$-integral on $%
\mathbb{Z}
_{p}$ is defined by T. Kim:%
\begin{equation}
I_{-q}\left( f\right) =\int_{%
\mathbb{Z}
_{p}}f\left( \xi \right) d\mu _{-q}\left( \xi \right) =\lim_{N\rightarrow
\infty }\sum_{\xi =0}^{p^{N}-1}q^{\xi }f\left( \xi \right) \left( -1\right)
^{\xi }  \label{equation 1}
\end{equation}

(for more informations on this subject, see \cite{kim 14}, \cite{kim 15} and 
\cite{kim 16}).

From (\ref{equation 1}), we have well known the following equality:%
\begin{equation}
qI_{-q}\left( f_{1}\right) +I_{-q}\left( f\right) =\left[ 2\right]
_{q}f\left( 0\right)  \label{equation 4}
\end{equation}

here $f_{1}\left( x\right) :=f\left( x+1\right) $ (for details, see[2-40]).

Let $C\left( \left[ 0,1\right] \right) $ be the space of continuous
functions on $\left[ 0,1\right] $. For $C\left( \left[ 0,1\right] \right) $,
the weighted $q$-Bernstein operator for $f$ is defined by 
\begin{equation*}
\mathcal{B}_{n,q}^{\left( \alpha \right) }\left( f,x\right)
=\sum_{k=0}^{n}f\left( \frac{k}{n}\right) B_{k,n}^{\left( \alpha \right)
}\left( x\mid q\right) =\sum_{k=0}^{n}f\left( \frac{k}{n}\right) \binom{n}{k}%
\left[ x\right] _{q^{\alpha }}^{k}\left[ 1-x\right] _{q^{-\alpha }}^{n-k}
\end{equation*}

where $n,$ $k\in 
\mathbb{N}
^{\ast }$. Here $B_{k,n}^{\left( \alpha \right) }\left( x\mid q\right) $ is
called weighted $q$-Bernstein polynomials, which are defined by 
\begin{equation}
B_{k,n}^{\left( \alpha \right) }\left( x\mid q\right) =\binom{n}{k}\left[ x%
\right] _{q^{\alpha }}^{k}\left[ 1-x\right] _{q^{-\alpha }}^{n-k}\text{, }%
x\in \left[ 0,1\right]   \label{equation 2}
\end{equation}

(for more informations on this subject, see \cite{Araci 1}, \cite{kim 19}, 
\cite{Ryoo 1} and \cite{Ryoo 2}).

As is well known, the ordinary Genocchi polynomials are defined by menas of
the following generating function:%
\begin{equation}
\sum_{n=0}^{\infty }G_{n}\left( x\right) \frac{t^{n}}{n!}=e^{G\left(
x\right) t}=\frac{2t}{e^{t}+1}e^{xt}\text{.}  \label{equation 19}
\end{equation}

where the usual convention about replacing $G^{n}\left( x\right) $ by $%
G_{n}\left( x\right) $. For $x=0$ in (\ref{equation 19}), we have to $%
G_{n}\left( 0\right) :=G_{n}$, which is called Genocchi numbers. Then, we
can write the following%
\begin{equation}
e^{Gt}=\sum_{n=0}^{\infty }G_{n}\frac{t^{n}}{n!}=\frac{2t}{e^{t}+1}\text{.}
\label{equation 20}
\end{equation}

In \cite{Araci 2}, the $q$-Genocchi numbers are given as%
\begin{equation*}
G_{0,q}=0\text{ and }q\left( qG_{q}+1\right) ^{n}+G_{n,q}=\left\{ 
\begin{array}{cc}
\left[ 2\right] _{q} & \text{if }n=1 \\ 
0 & \text{if }n\neq 1%
\end{array}%
\right.
\end{equation*}

where the usual convention about replacing $\left( G_{q}\right) ^{n}$ by $%
G_{n,q}$.

For any $n\in 
\mathbb{N}
^{\ast }$, the ($h,q$)-Genocchi numbers are defined by 
\begin{equation*}
G_{0,q}^{\left( h\right) }=0\text{ and }q^{h-1}\left( qG_{q}^{\left(
h\right) }+1\right) ^{n}+G_{n,q}^{\left( h\right) }=\left\{ 
\begin{array}{cc}
\left[ 2\right] _{q} & \text{if }n=1 \\ 
0 & \text{if }n\neq 1%
\end{array}%
\right.
\end{equation*}

where the usual convention about replacing $\left( G_{q}^{\left( h\right)
}\right) ^{n}$ by $G_{n,q}^{\left( h\right) }$ (for details, see \cite{Araci
10}).

Recently, Araci $et$ $al$. are defined the $(h,q)$-Genocchi numbers with
weight $\alpha $ by 
\begin{equation}
\frac{\widetilde{G}_{n+1,q}^{\left( \alpha ,h\right) }\left( x\right) }{n+1}%
=\int_{%
\mathbb{Z}
_{p}}q^{\left( h-1\right) \xi }\left[ x+\xi \right] _{q^{\alpha }}^{n}d\mu
_{-q}\left( \xi \right) \text{.}  \label{equation 21}
\end{equation}

By (\ref{equation 21}), we have the following identity%
\begin{equation}
\widetilde{G}_{n,q}^{\left( \alpha ,h\right) }\left( x\right) =\sum_{k=0}^{n}%
\binom{n}{k}q^{\alpha kx}\widetilde{G}_{n,q}^{\left( \alpha ,h\right) }\left[
x\right] _{q^{\alpha }}^{n-k}=q^{-\alpha x}\left( q^{\alpha x}\widetilde{G}%
_{q}^{\left( \alpha ,h\right) }+\left[ x\right] _{q^{\alpha }}\right) ^{n}
\label{equation 23}
\end{equation}

where the usual convention about replacing $\left( \widetilde{G}_{q}^{\left(
\alpha ,h\right) }\right) ^{n}$ by $\widetilde{G}_{n,q}^{\left( \alpha
,h\right) }$ is used (for details, \cite{Araci 4}).

In this paper, we derive some new properties ($h,q$)-Genocchi numbers and
polynomials from the fermionic $p$-adic $q$-integral on $%
\mathbb{Z}
_{p}$. Also, we show that these type polynomials are related to ($h,q$%
)-Genocchi numbers and polynomials.

\section{\textbf{On the }$\left( h,q\right) $-\textbf{Genocchi numbers\ and
polynomials}}

In this section, we consider the ($h,q$)-Genocchi numbers and polynomials by
using fermionic $p$-adic $q$-integral on $%
\mathbb{Z}
_{p}$ and the weighted $q$-Bernstein polynomials. We can now start the
following expression.

In \cite{Araci 4}, we have the ($h,q$)-Genocchi numbers as follows: For $%
\alpha \in 
\mathbb{N}
^{\ast }$ and $n,h\in 
\mathbb{N}
$,%
\begin{equation}
\widetilde{G}_{0,q}^{\left( \alpha ,h\right) }=0\text{ and }q^{h}\widetilde{G%
}_{n,q}^{\left( \alpha ,h\right) }\left( 1\right) +\widetilde{G}%
_{n,q}^{\left( \alpha ,h\right) }=\left\{ 
\begin{array}{cc}
\left[ 2\right] _{q} & \text{if }n=1, \\ 
0 & \text{if }n\neq 1.%
\end{array}%
\right.  \label{equation 22}
\end{equation}

By (\ref{equation 23}) and (\ref{equation 22}), we have the following
corollary.

\begin{corollary}
For $\alpha \in 
\mathbb{N}
^{\ast }$ and $n,h\in 
\mathbb{N}
$, then we have%
\begin{equation}
\widetilde{G}_{0,q}^{\left( \alpha ,h\right) }=0\text{ and }q^{h-\alpha
}\left( q^{\alpha }\widetilde{G}_{q}^{\left( \alpha ,h\right) }+1\right)
^{n}+\widetilde{G}_{n,q}^{\left( \alpha ,h\right) }=\left\{ 
\begin{array}{cc}
\left[ 2\right] _{q} & \text{if }n=1, \\ 
0 & \text{if }n\neq 1.%
\end{array}%
\right.  \label{equation 25}
\end{equation}
\end{corollary}

By (\ref{equation 21}), we get symmetric property that%
\begin{eqnarray*}
\frac{\widetilde{G}_{n+1,q^{-1}}^{\left( \alpha ,h\right) }\left( 1-x\right) 
}{n+1} &=&\int_{%
\mathbb{Z}
_{p}}q^{\left( 1-h\right) \xi }\left[ 1-x+\xi \right] _{q^{-\alpha
}}^{n}d\mu _{-q^{-1}}\left( \xi \right) \\
&=&\left( -1\right) ^{n}q^{h+\alpha n-1}\int_{%
\mathbb{Z}
_{p}}q^{\left( h-1\right) \xi }\left[ x+\xi \right] _{q^{\alpha }}^{n}d\mu
_{-q}\left( \xi \right)
\end{eqnarray*}

From this, we state the following theorem.

\begin{theorem}
The following identity%
\begin{equation}
\widetilde{G}_{n+1,q^{-1}}^{\left( \alpha ,h\right) }\left( 1-x\right)
=\left( -1\right) ^{n}q^{h+\alpha n-1}\widetilde{G}_{n+1,q}^{\left( \alpha
,h\right) }\left( x\right)  \label{equation 24}
\end{equation}%
is true.
\end{theorem}

By using (\ref{equation 23}), (\ref{equation 22}) and (\ref{equation 25}),
we compute as follows:%
\begin{eqnarray}
q^{2\alpha }\widetilde{G}_{n,q}^{\left( \alpha ,h\right) }\left( 2\right)
&=&\left( q^{2\alpha }\widetilde{G}_{q}^{\left( \alpha ,h\right) }+\left[ 2%
\right] _{q^{\alpha }}\right) ^{n}  \label{equation 6} \\
&=&\sum_{l=0}^{n}\binom{n}{l}q^{\alpha l}\left( q^{\alpha }\widetilde{G}%
_{q}^{\left( \alpha ,h\right) }+1\right) ^{l}  \notag \\
&=&nq^{2\alpha -h}\left( \left[ 2\right] _{q}-\widetilde{G}_{1,q}^{\left(
\alpha ,h\right) }\right) -q^{\alpha -h}\sum_{l=2}^{n}\binom{n}{l}q^{\alpha
l}\widetilde{G}_{l,q}^{\left( \alpha ,h\right) }  \notag \\
&=&nq^{2\alpha -h}\left[ 2\right] _{q}+q^{2\alpha -2h}\widetilde{G}%
_{n,q}^{\left( \alpha ,h\right) }\text{ if }n>1.  \notag
\end{eqnarray}

After the above applications, we procure the following theorem.

\begin{theorem}
For $n>1$, then we have%
\begin{equation*}
\widetilde{G}_{n,q}^{\left( \alpha ,h\right) }\left( 2\right) =nq^{-h}\left[
2\right] _{q}+q^{-2h}\widetilde{G}_{n,q}^{\left( \alpha ,h\right) }.
\end{equation*}
\end{theorem}

We need the following equality for sequel of this paper:%
\begin{equation}
\left[ 1-x\right] _{q^{-\alpha }}^{n}=\left( \frac{1-q^{-\alpha \left(
1-x\right) }}{1-q^{-\alpha }}\right) ^{n}=\left( -1\right) ^{n}q^{n\alpha }%
\left[ x-1\right] _{q^{\alpha }}^{n}\text{.}  \label{equation 26}
\end{equation}

Now also, by using (\ref{equation 26}), we consider the following%
\begin{eqnarray*}
&&q^{h-1}\int_{%
\mathbb{Z}
_{p}}q^{\left( h-1\right) \xi }\left[ 1-\xi \right] _{q^{-\alpha }}^{n}d\mu
_{-q}\left( \xi \right) \\
&=&\left( -1\right) ^{n}q^{h+n\alpha -1}\int_{%
\mathbb{Z}
_{p}}q^{\left( h-1\right) \xi }\left[ \xi -1\right] _{q^{\alpha }}^{n}d\mu
_{-q}\left( \xi \right) \\
&=&\left( -1\right) ^{n}q^{h+n\alpha -1}\frac{\widetilde{G}_{n+1,q}^{\left(
\alpha ,h\right) }\left( -1\right) }{n+1}\text{.}
\end{eqnarray*}

By considering last identity and (\ref{equation 24}), we get the following
theorem.

\begin{theorem}
The following identity holds true:%
\begin{equation}
\int_{%
\mathbb{Z}
_{p}}q^{\left( h-1\right) \left( \xi +1\right) }\left[ 1-\xi \right]
_{q^{-\alpha }}^{n}d\mu _{-q}\left( \xi \right) =\frac{\widetilde{G}%
_{n+1,q^{-1}}^{\left( \alpha ,h\right) }\left( 2\right) }{n+1}\text{.}
\label{equation 27}
\end{equation}
\end{theorem}

From (\ref{equation 27}), we have the following%
\begin{equation}
\int_{%
\mathbb{Z}
_{p}}q^{\left( h-1\right) \xi }\left[ 1-\xi \right] _{q^{-\alpha }}^{n}d\mu
_{-q}\left( \xi \right) =\left[ 2\right] _{q}+q^{h+1}\frac{\widetilde{G}%
_{n+1,q^{-1}}^{\left( \alpha ,h\right) }}{n+1}\text{.}  \notag
\end{equation}

Thus, we obtain the following theorem.

\begin{theorem}
The following identity 
\begin{equation}
\int_{%
\mathbb{Z}
_{p}}q^{\left( h-1\right) \xi }\left[ 1-\xi \right] _{q^{-\alpha }}^{n}d\mu
_{-q}\left( \xi \right) =\left[ 2\right] _{q}+q^{h+1}\frac{\widetilde{G}%
_{n+1,q^{-1}}^{\left( \alpha ,h\right) }}{n+1}  \label{equation 15}
\end{equation}%
is true.
\end{theorem}

\section{\textbf{Some new identities on the }$\left( h,q\right) $\textbf{%
-Genocchi numbers}}

\qquad In this section, we introduce new identities of the ($h,q$)-Genocchi
numbers, that is, we derive some interesting and worthwhile relations for
studying in Theory of Analytic Numbers.

\qquad For $x\in \left[ 0,1\right] $, we give definition of weighted $q$%
-Bernstein polynomials as follows:%
\begin{equation}
B_{k,n}^{\left( \alpha \right) }\left( x\mid q\right) =\binom{n}{k}\left[ x%
\right] _{q^{\alpha }}^{k}\left[ 1-x\right] _{q^{-\alpha }}^{n-k}\text{,
where }n,k\in 
\mathbb{Z}
_{+}\text{.}  \label{equation 9}
\end{equation}

By expression of (\ref{equation 9}), we have the properties of symmetry of
weighted $q$-Bernstein polynomials as follows:%
\begin{equation}
B_{k,n}^{\left( \alpha \right) }\left( x\mid q\right) =B_{n-k,n}^{\left(
\alpha \right) }\left( 1-x\mid \frac{1}{q}\right) \text{, (for details, see 
\cite{kim 19}).}  \label{equation 10}
\end{equation}

Thus, (\ref{equation 15}), (\ref{equation 9}) and (\ref{equation 10}), we
see that%
\begin{eqnarray*}
I_{1} &=&\int_{%
\mathbb{Z}
_{p}}q^{\left( h-1\right) x}B_{k,n}^{\left( \alpha \right) }\left( x\mid
q\right) d\mu _{-q}\left( x\right) =\binom{n}{k}\int_{%
\mathbb{Z}
_{p}}q^{\left( h-1\right) x}\left[ x\right] _{q^{\alpha }}^{k}\left[ 1-x%
\right] _{q^{-\alpha }}^{n-k}d\mu _{-q}\left( x\right) \\
&=&\binom{n}{k}\sum_{l=0}^{k}\binom{k}{l}\left( -1\right) ^{k+l}\int_{%
\mathbb{Z}
_{p}}q^{\left( h-1\right) x}\left[ 1-x\right] _{q^{-\alpha }}^{n-l}d\mu
_{-q}\left( x\right) \\
&=&\binom{n}{k}\sum_{l=0}^{k}\binom{k}{l}\left( -1\right) ^{k+l}\left\{ %
\left[ 2\right] _{q}+q^{h+1}\frac{\widetilde{G}_{n-l+1,q^{-1}}^{\left(
\alpha ,h\right) }}{n-l+1}\right\} \\
&=&\left\{ 
\begin{array}{cc}
\left[ 2\right] _{q}+q^{h+1}\frac{\widetilde{G}_{n+1,q^{-1}}^{\left( \alpha
,h\right) }}{n+1} & \text{if }k=0, \\ 
\binom{n}{k}\sum_{l=0}^{k}\binom{k}{l}\left( -1\right) ^{k+l}\left\{ \left[ 2%
\right] _{q}+q^{h+1}\frac{\widetilde{G}_{n-l+1,q^{-1}}^{\left( \alpha
,h\right) }}{n-l+1}\right\} & \text{if }k\neq 0.%
\end{array}%
\right.
\end{eqnarray*}

On the other hand, for $n$, $k\in 
\mathbb{Z}
_{+}$ with $n>k$, we compute%
\begin{eqnarray*}
I_{2} &=&\int_{%
\mathbb{Z}
_{p}}q^{\left( h-1\right) x}B_{k,n}^{\left( \alpha \right) }\left( x\mid
q\right) d\mu _{-q}\left( x\right) \\
&=&\binom{n}{k}\int_{%
\mathbb{Z}
_{p}}q^{\left( h-1\right) x}\left[ x\right] _{q^{\alpha }}^{k}\left[ 1-x%
\right] _{q^{-\alpha }}^{n-k}d\mu _{-q}\left( x\right) \\
&=&\binom{n}{k}\sum_{l=0}^{n-k}\binom{n-k}{l}\left( -1\right) ^{l}\int_{%
\mathbb{Z}
_{p}}q^{\left( h-1\right) x}\left[ x\right] _{q^{\alpha }}^{l+k}d\mu
_{-q}\left( x\right) \\
&=&\binom{n}{k}\sum_{l=0}^{n-k}\binom{n-k}{l}\left( -1\right) ^{l}\frac{%
\widetilde{G}_{l+k+1,q}^{\left( \alpha ,h\right) }}{l+k+1}\text{.}
\end{eqnarray*}

Equating $I_{1}$ and $I_{2}$, then we have the following theorem.

\begin{theorem}
The following identity holds true:%
\begin{equation*}
\sum_{l=0}^{n-k}\binom{n-k}{l}\left( -1\right) ^{l}\frac{\widetilde{G}%
_{l+k+1,q}^{\left( \alpha ,h\right) }}{l+k+1}=\left\{ 
\begin{array}{cc}
\left[ 2\right] _{q}+q^{h+1}\frac{\widetilde{G}_{n+1,q^{-1}}^{\left( \alpha
,h\right) }}{n+1} & \text{if }k=0, \\ 
\sum_{l=0}^{k}\binom{k}{l}\left( -1\right) ^{k+l}\left\{ \left[ 2\right]
_{q}+q^{h+1}\frac{\widetilde{G}_{n-l+1,q^{-1}}^{\left( \alpha ,h\right) }}{%
n-l+1}\right\} & \text{if }k\neq 0.%
\end{array}%
\right.
\end{equation*}
\end{theorem}

Let $n_{1},n_{2},k\in 
\mathbb{Z}
_{+}$ with $n_{1}+n_{2}>2k$. Then, we derive the followings%
\begin{eqnarray}
I_{3} &=&\int_{%
\mathbb{Z}
_{p}}q^{\left( h-1\right) x}B_{k,n_{1}}^{\left( \alpha \right) }\left( x\mid
q\right) B_{k,n_{2}}^{\left( \alpha \right) }\left( x\mid q\right) d\mu
_{-q}\left( x\right)  \notag \\
&=&\binom{n_{1}}{k}\binom{n_{2}}{k}\sum_{l=0}^{2k}\binom{2k}{l}\left(
-1\right) ^{2k+l}\int_{%
\mathbb{Z}
_{p}}q^{\left( h-1\right) x}\left[ 1-x\right] _{q^{-\alpha
}}^{n_{1}+n_{2}-l}d\mu _{-q}\left( x\right)  \notag \\
&=&\left( \binom{n_{1}}{k}\binom{n_{2}}{k}\sum_{l=0}^{2k}\binom{2k}{l}\left(
-1\right) ^{2k+l}\left( \left[ 2\right] _{q}+q^{h+1}\frac{\widetilde{G}%
_{n_{1}+n_{2}-l+1,q^{-1}}^{\left( \alpha ,h\right) }}{n_{1}+n_{2}-l+1}%
\right) \right)  \notag \\
&=&\left\{ 
\begin{array}{cc}
\left[ 2\right] _{q}+q^{h+1}\frac{\widetilde{G}_{n_{1}+n_{2}+1,q^{-1}}^{%
\left( \alpha ,h\right) }}{n+1} & \text{if }k=0, \\ 
\binom{n}{k}\sum_{l=0}^{2k}\binom{2k}{l}\left( -1\right) ^{2k+l}\left\{ %
\left[ 2\right] _{q}+q^{h+1}\frac{\widetilde{G}_{n_{1}+n_{2}-l+1,q^{-1}}^{%
\left( \alpha ,h\right) }}{n_{1}+n_{2}-l+1}\right\} & \text{if }k\neq 0.%
\end{array}%
\right.  \notag
\end{eqnarray}

In other words, by using the binomial theorem, we can derive the following
equation.%
\begin{eqnarray*}
I_{4} &=&\int_{%
\mathbb{Z}
_{p}}q^{\left( h-1\right) x}B_{k,n_{1}}^{\left( \alpha \right) }\left( x\mid
q\right) B_{k,n_{2}}^{\left( \alpha \right) }\left( x\mid q\right) d\mu
_{-q}\left( x\right) \\
&=&\dprod\limits_{i=1}^{2}\binom{n_{i}}{k}\sum_{l=0}^{n_{1}+n_{2}-2k}\binom{%
n_{1}+n_{2}-2k}{l}\left( -1\right) ^{l}\int_{%
\mathbb{Z}
_{p}}q^{\left( h-1\right) x}\left[ x\right] _{q^{\alpha }}^{2k+l}d\mu
_{-q}\left( x\right) \\
&=&\dprod\limits_{i=1}^{2}\binom{n_{i}}{k}\sum_{l=0}^{n_{1}+n_{2}-2k}\binom{%
n_{1}+n_{2}-2k}{l}\left( -1\right) ^{l}\frac{\widetilde{G}%
_{l+2k+1,q}^{\left( \alpha ,h\right) }}{l+2k+1}.
\end{eqnarray*}

Combining $I_{3}$ and $I_{4}$, we state the following theorem.

\begin{theorem}
For $n_{1},n_{2},k\in 
\mathbb{Z}
_{+}$ with $n_{1}+n_{2}>2k,$ we have%
\begin{eqnarray*}
&&\sum_{l=0}^{n_{1}+n_{2}-2k}\binom{n_{1}+n_{2}-2k}{l}\left( -1\right) ^{l}%
\frac{\widetilde{G}_{l+2k+1,q}^{\left( \alpha ,h\right) }}{l+2k+1} \\
&=&\left\{ 
\begin{array}{cc}
\left[ 2\right] _{q}+q^{h+1}\frac{\widetilde{G}_{n_{1}+n_{2}+1,q^{-1}}^{%
\left( \alpha ,h\right) }}{n_{1}+n_{2}+1} & \text{if }k=0, \\ 
\sum_{l=0}^{2k}\binom{2k}{l}\left( -1\right) ^{2k+l}\left\{ \left[ 2\right]
_{q}+q^{h+1}\frac{\widetilde{G}_{n_{1}+n_{2}-l+1,q^{-1}}^{\left( \alpha
,h\right) }}{n_{1}+n_{2}-l+1}\right\} & \text{if }k\neq 0.%
\end{array}%
\right.
\end{eqnarray*}
\end{theorem}

For $x\in 
\mathbb{Z}
_{p}$ and $s\in 
\mathbb{N}
$ with $s\geq 2,$ let $n_{1},n_{2},...,n_{s},k\in 
\mathbb{Z}
_{+}$ with $\sum_{l=1}^{s}n_{l}>sk$. Then we take the fermionic $p$-adic $q$%
-integral on $%
\mathbb{Z}
_{p}$ for the weighted $q$-Bernstein polynomials of degree $n$ as follows:%
\begin{eqnarray*}
I_{5} &=&\int_{%
\mathbb{Z}
_{p}}q^{\left( h-1\right) x}\left\{
\tprod\limits_{i=1}^{s}B_{k,n_{i}}^{\left( \alpha \right) }\left( x\mid
q\right) \right\} d\mu _{-q}\left( x\right) \\
&=&\dprod\limits_{i=1}^{s}\binom{n_{i}}{k}\int_{%
\mathbb{Z}
_{p}}\left[ x\right] _{q^{\alpha }}^{sk}\left[ 1-x\right] _{q^{-\alpha
}}^{n_{1}+n_{2}+...+n_{s}-sk}q^{\left( h-1\right) x}d\mu _{-q}\left( x\right)
\\
&=&\dprod\limits_{i=1}^{s}\binom{n_{i}}{k}\sum_{l=0}^{sk}\binom{sk}{l}\left(
-1\right) ^{l+sk}\int_{%
\mathbb{Z}
_{p}}\left[ 1-x\right] _{q^{-\alpha }}^{n_{1}+n_{2}+...+n_{s}-l}q^{\left(
h-1\right) x}d\mu _{-q}\left( x\right) \\
&=&\left\{ 
\begin{array}{cc}
\left[ 2\right] _{q}+q^{h+1}\frac{\widetilde{G}%
_{n_{1}+n_{2}+...+n_{s}+1,q^{-1}}^{\left( \alpha ,h\right) }}{%
n_{1}+n_{2}+...+n_{s}+1} & \text{if }k=0, \\ 
\dprod\limits_{i=1}^{s}\binom{n_{i}}{k}\sum_{l=0}^{sk}\binom{sk}{l}\left(
-1\right) ^{sk+l}\left\{ \left[ 2\right] _{q}+q^{h+1}\frac{\widetilde{G}%
_{n_{1}+n_{2}+...+n_{s}-l+1,q^{-1}}^{\left( \alpha ,h\right) }}{%
n_{1}+n_{2}+...+n_{s}-l+1}\right\} & \text{if }k\neq 0.%
\end{array}%
\right.
\end{eqnarray*}

On the other hand, from the definition of weighted $q$-Bernstein polynomials
and the binomial theorem, we easily get%
\begin{eqnarray*}
I_{6} &=&\int_{%
\mathbb{Z}
_{p}}q^{\left( h-1\right) x}\left\{
\tprod\limits_{i=1}^{s}B_{k,n_{i}}^{\left( \alpha \right) }\left( x\mid
q\right) \right\} d\mu _{-q}\left( x\right) \\
&=&\dprod\limits_{i=1}^{s}\binom{n_{i}}{k}\sum_{l=0}^{n_{1}+...+n_{s}-sk}%
\binom{\sum_{d=1}^{s}\left( n_{d}-k\right) }{l}\left( -1\right) ^{l}\int_{%
\mathbb{Z}
_{p}}\left[ x\right] _{q^{\alpha }}^{sk+l}q^{\left( h-1\right) x}d\mu
_{-q}\left( x\right) \\
&=&\dprod\limits_{i=1}^{s}\binom{n_{i}}{k}\sum_{l=0}^{n_{1}+...+n_{s}-sk}%
\binom{\sum_{d=1}^{s}\left( n_{d}-k\right) }{l}\left( -1\right) ^{l}\frac{%
\widetilde{G}_{l+sk+1,q}^{\left( \alpha ,h\right) }}{l+sk+1}\text{.}
\end{eqnarray*}

Equating $I_{5}$ and $I_{6}$, we discover the following theorem.

\begin{theorem}
For $s\in 
\mathbb{N}
$ with $s\geq 2$, let $n_{1},n_{2},...,n_{s},k\in 
\mathbb{Z}
_{+}$ with $\sum_{l=1}^{s}n_{l}>sk.$ Then, we have 
\begin{eqnarray*}
&&\sum_{l=0}^{n_{1}+...+n_{s}-sk}\binom{\sum_{d=1}^{s}\left( n_{d}-k\right) 
}{l}\left( -1\right) ^{l}\frac{\widetilde{G}_{l+sk+1,q}^{\left( \alpha
,h\right) }}{l+sk+1} \\
&=&\left\{ 
\begin{array}{cc}
\left[ 2\right] _{q}+q^{h+1}\frac{\widetilde{G}%
_{n_{1}+n_{2}+...+n_{s}+1,q^{-1}}^{\left( \alpha ,h\right) }}{%
n_{1}+n_{2}+...+n_{s}+1} & \text{if }k=0, \\ 
\sum_{l=0}^{sk}\binom{sk}{l}\left( -1\right) ^{sk+l}\left\{ \left[ 2\right]
_{q}+q^{h+1}\frac{\widetilde{G}_{n_{1}+n_{2}+...+n_{s}-l+1,q^{-1}}^{\left(
\alpha ,h\right) }}{n_{1}+n_{2}+...+n_{s}-l+1}\right\} & \text{if }k\neq 0.%
\end{array}%
\right.
\end{eqnarray*}
\end{theorem}

\end{document}